\documentclass{smjour}
 \usepackage{combdesi} 

\usepackage[T1]{fontenc}

\newcommand{\D}{\displaystyle}

\newcommand{\Z}{\mathsf{\textstyle Z\kern-0.4em Z}}

\begin{document}



\yearofpublication{(Year)}
\volume{(Volume Number)}
\issuenumber{(Optional Issue Number)}
\cccline{(cccline information)}


\authorrunninghead{Alvarez, Gudiel and Guemes}
\titlerunninghead{On $\Z _{\lowercase{t}} \times \Z_2^2$-cocyclic Hadamard matrices}

\setcounter{page}{1} 



\title{On $\Z _t \times \Z_2^2$-cocyclic Hadamard matrices}

\author{V\'{\i}ctor~\'Alvarez,
                }

\affil{Department of Applied Mathematics 1, University of Seville,
Seville, Spain e-mail: valvarez@us.es}
\author{F\'{e}lix~Gudiel,
        }

\affil{Department of Applied Mathematics 1, University of Seville,
Seville, Spain e-mail: gudiel@us.es}

\author{Maria~Bel\'{e}n~G\"{u}emes }
\affil{Department of Algebra, University of Seville,
Seville, Spain e-mail: bguemes@us.es}


\abstract{A characterization of $\Z _t \times \Z_2^2$-cocyclic
Hadamard matrices is described, depending on the notions of {\em
distributions}, {\em ingredients} and {\em recipes}. In
particular, these notions lead to the establishment of some bounds
on the number and distribution of 2-coboundaries over $\Z_t \times
\Z _2^2$ to use and the way in which they have to be combined in
order to obtain a $\Z _t \times \Z_2^2$-cocyclic Hadamard matrix.
Exhaustive searches have been performed, so that the table in p.
132 in \cite{BH95} is corrected and completed. Furthermore, we
identify four different operations  on the set of coboundaries
defining $\Z _t \times \Z_2^2$-cocyclic matrices, which  preserve
orthogonality. We split the set of
Hadamard matrices into disjoint orbits, define representatives for
them and take advantage of this fact to compute them in an easier
way than the usual purely exhaustive way, in terms of {\em diagrams}.
Let ${\cal H}$ be the set of cocyclic Hadamard matrices over $\Z_t
\times \Z_2^2$ having a symmetric diagram. We also prove that the set of
Williamson type matrices is a  subset of
${\cal H}$ of size $\frac{|{\cal H}|}{t}$.}


\begin{article} 


\section{Introduction}

Hadamard matrices are $n\times n$ square matrices $H$ with entries
in $\{1, -1\}$ such that every pair of rows (respectively,
columns) are orthogonal, that is, $H  H^T=nI_n$.

Due to this nice combinatorial property, Hadamard matrices have
many applications in a wide variety of fields, such as Signal
Processing, Coding Theory and Cryptography (see \cite{Hor07} for
details). Consequently, there is a real interest in knowing enough
Hadamard matrices for practical use.

It is a straightforward exercise to prove that the order of a
Hadamard matrix has to be 1, 2 or a multiple of 4 (as soon as
three or more rows have to be simultaneously orthogonal one to
each other). Unfortunately, the Hadamard Conjecture about the
existence of these matrices for every order $4t$ remains unproved
since the XIXth Century.

Nowadays, there are three orders less than 1000 for which no
Hadamard matrix is known: $668=4 \cdot 167$, $716=4\cdot 179$, and
$892=4 \cdot 223$. Furthermore, there are 9 orders in the range
$[1000,2000]$ for which no Hadamard matrix is known (see
\cite{Hor07,DGK13} for details).

One of the most promising techniques for constructing Hadamard
matrices is the cocyclic approach (see \cite{HdL93,HdL95,Hor07}).
A cocyclic matrix $M_f$ over a group $G$ is a matrix
$M_f=(f(g,h))$, for $f$ a $2$-cocycle over $G$, that is, a
function $f:G \times G \rightarrow \{1,-1\}$ such that for every
$a,b,c \in G$, $f(a,b)\cdot f(ab,c)\cdot f(a,bc) \cdot f(b,c)=1$.

Actually, many well known families of Hadamard matrices, such as
Sylvester's, Paley's, Williamson's and Ito's, have shown to be
cocyclic over appropriate groups (see \cite{Hor07} for details).
This has provided inspiration for the Cocyclic Hadamard
Conjecture, which states that cocyclic Hadamard matrices exist for
every order $4t$.

In this paper we are interested in characterizing cocyclic
Hadamard matrices over $\Z_t \times \Z_2^2$, which include the
family of symmetric Williamson type Hadamard matrices.

Following the indications of \cite{AAFR08}, we will describe
bounds on the number of 2-coboundaries over $\Z_t \times \Z_2^2$
to be combined, as well as their distribution (in terms of what we
call {\em ingredients} and {\em recipes}), in order to construct a
$\Z_t \times \Z_2^2$-cocyclic Hadamard matrix.

This information will   allow us to design an exhaustive search
for $\Z_t \times \Z_2^2$-cocyclic Hadamard matrices for $3 \leq t
\leq 13$, so that the table in p. 132 in \cite{BH95} is corrected
and completed.

Next, we will introduce what we call {\it diagrams}, a visual
representation of the coboundaries which define a $\Z _t \times
\Z_2^2$-cocyclic matrix and  we will study four different
operations  on the set of coboundary matrices over $\Z_t \times
\Z_2^2$: complements, rotations, swappings and dilatations. In
particular, these operations extend to operations over $\Z_t
\times \Z_2^2$-cocyclic matrices, which will be proved to preserve
orthogonality. These operations partition
the set of $\Z_t \times \Z_2^2$-cocyclic matrices into disjoint
orbits, which can be easily computed once one element is known.
Among all the elements of an orbit, a representative can be
chosen, in a standard way that will be made precise.

Finally, by applying these ideas to the task of searching for
cocyclic Hadamard matrices over $\Z _t \times \Z_2^2$, we have
been able to extend the table in p. 132 in \cite{BH95} for values
of $t$ in the range $3 \leq t \leq 23$, and to explain the fact
that the set of symmetric Williamson type Hadamard matrices
obtained from symmetric diagrams is in proportion $\frac{1}{t}$
with respect to the full set of  $\Z _t \times \Z_2^2$-cocyclic
Hadamard matrices.

We organize the paper as follows. Section 2 
is dedicated to
describe all about $\Z_t \times \Z_2^2$-cocyclic Hadamard
matrices. In Section 3 
we introduce the notions of {\em
distribution}, {\em ingredients} and {\em recipes}, in terms of
which we find some upper and lower bounds on the number of
2-coboundaries over $\Z_t \times \Z_2^2$ which have to be combined
in order to get $\Z_t \times \Z_2^2$-cocyclic Hadamard matrices,
as well as the way in which they have to be distributed, and the
results obtained. Section 4 
is dedicated to the description
of diagrams, after a discussion about the convenience of using the
whole set of coboundaries instead of the basis, in order to
represent $\Z _t \times \Z_2^2$-cocyclic matrices. Section 5 
defines four operations on $\Z_t \times \Z_2^2$-cocyclic matrices which preserve orthogonality, and which have a nice
interpretation in terms of diagrams. These operations split the
set of $\Z_t \times \Z_2^2$-cocyclic Hadamard matrices into disjoint orbits, which can be
generated from any of their elements (for instance by its
representative). In Section 6 
we focus on the case of $\Z_t \times \Z_2^2$-cocyclic Hadamard matrices obtained from symmetric diagrams, which in turn permit counting the number of Williamson type Hadamard matrices. We also include some final remarks and
further work.

\section{$\Z_{\lowercase{t}} \times \Z_2^2$-cocyclic Hadamard matrices}\label{2}

Consider the group $G=\Z_t \times \Z_2^2 $, a  basis ${\cal
B}=\{\partial _2,\ldots, \partial
_{4t-2},\beta_1,\beta_2,\gamma\}$ for 2-cocycles over $G$ is
described in \cite{AAFR09}, and consists of $4t-3$ coboundaries
$\partial _k$, two cocycles $\beta_i$ coming from inflation and
one cocycle $\gamma$ coming from transgression.


It has been observed that cocyclic Hadamard matrices over $\Z _t
\times \Z_2^2$ mostly use all the three representative cocycles
$\beta_1,\beta_2$ and $\gamma$ simultaneously (see \cite{BH95} for
details). We will assume that every cocyclic matrix $M$ is
obtained as a product $M=M_{\partial_{i_1}} \ldots
M_{\partial_{i_w}}\cdot R$,  for  $2 \leq i_1 < \ldots <i_w \leq
4t-2$, where $R=M_{\beta_1} \cdot M_{\beta_2} \cdot M_{\gamma}$.

Here $M_{\partial _i}$ refers to the {\em generalized} coboundary matrix
associated to the $i^{th}$-element in $G$, with the $i^{th}$-row
and the $i^{th}$-column negated, as introduced in  \cite{AAFR09}.

%

In particular, there are three coboundary matrices which are not
in ${\cal B}$: $M_{\partial _1}$, $M_{\partial _{4t-1}}$ and
$M_{\partial _{4t}}$.
%
%
%
%
%
%
%
Consequently, every $\Z _t \times \Z_2^2$-cocyclic matrix $M_f$
using $R$ may be expressed as a pointwise product of matrices in
$\{ M_{\partial _1},\ldots ,M_{\partial _{4t}}\}$ in 8 different
ways,  just one of which does not use any  of $M_{\partial _1}$,
$M_{\partial _{4t-1}}$ and $M_{\partial _{4t}}$ (and gives
precisely the expression of $M_f$ as a linear combination of
elements in ${\cal B}$). Actually, suppose that $M_f= \D R\cdot
\prod _{k=1}^4 \prod _{i_j \in J_k} M_{\partial_{i_j}}$, where
$J_k \subset \{1, \ldots, 4t\}$ is a subset of indexes which are
congruent to $k$ modulo 4. Then $M_f$ may be expressed as the
pointwise product of the coboundary matrices of indexes belonging
to any of the following 8 subsets:
$(J_1,J_2,J_3,J_4)$, $(\bar{J_1},\bar{J_2},J_3,J_4)$,
$(J_1,\bar{J_2},\bar{J_3},J_4)$, $(J_1,\bar{J_2},J_3,\bar{J_4})$,
$(\bar{J_1},J_2,\bar{J_3},J_4)$,
$(\bar{J_1},J_2,J_3,\bar{J_4})$, $(J_1,J_2,\bar{J_3},\bar{J_4})$,
$(\bar{J_1},\bar{J_2},\bar{J_3},\bar{J_4})$, where
$\bar{J_k}=\{k+4i: 0 \leq i \leq t-1\} \backslash J_k$ denotes the
complementary subset of $J_k$.

It is known that a cocyclic matrix is Hadamard if and only if the
summation of each row but the first is zero (this is the cocyclic
Hadamard test, see \cite{HdL95,BH95} for instance). Furthermore,
as proved in \cite{AAFR08}, a $\Z _t \times \Z_2^2$-cocyclic
matrix is Hadamard if and only if the summation of each of the
rows from 5 to $2t+2$ is 0, and an equivalent characterization of
the cocyclic Hadamard test may be described in terms of {\em n-paths}
($c_n$) and {\em n-intersections} ($I_n$) (the interested reader is referred to \cite{AAFR08} for a precise definition about paths and intersections).
In particular, the following result derives straightforwardly from the work in \cite{AAFR08}.

\begin{proposition} \label{corolario1} For  $\Z _t \times \Z_2^2$-cocyclic matrices:

\begin{enumerate}

\item the summation of a row $n \equiv 1 \ mod \ 4$ is zero if and only if \begin{equation}\label{filahadamard1} c_n=t, \end{equation}

\item the summation of a row $n \equiv 0,2,3 \ mod \ 4$ is zero if and only if \begin{equation}\label{filahadamard2} c_n=I_n. \end{equation}

\end{enumerate}
\end{proposition}


In the following section we analyze the number $c_n$ of $n$-paths
of a $\Z _t \times \Z_2^2$-cocyclic matrix. Focusing in rows $n
\equiv 1 \ mod \ 4$, we will obtain upper and lower bounds on the
number of coboundaries to combine in order to get a $\Z _t \times
\Z_2^2$-cocyclic Hadamard matrix. Furthermore, we will
characterize the distribution of these coboundaries in terms of
{\em ingredients} and {\em recipes}.

\section{Distributions, ingredients and recipes}\label{3}

Firstly, we analyze the way in which $n$-paths are generated on
$\Z _t \times \Z_2^2$-cocyclic matrices, depending on the value of
$n$ modulo $4$.

\begin{lemma} \label{lema1} Characterization of $n$-paths of coboundaries on $\Z _t \times \Z_2^2$-cocyclic matrices:
\begin{enumerate}
\item If $n \equiv 1 \ mod \ 4$,  $M_{\partial _{i-n+1}}$  forms an $n$-path with $M_{\partial _i}$.

\item If $n \equiv 2 \ mod \ 4$,  $M_{\partial _{i-n+2-(-1)^i}}$  forms an $n$-path with $M_{\partial _i}$.

\item If $n \equiv 3 \ mod \ 4$, $M_{\partial _{i-n+3-2(-1)^{\lceil \frac{i \ mod \ 4}{2} \rceil}}}$ forms an $n$-path with $M_{\partial _i}$.

\item If $n \equiv 0 \ mod \ 4$, $M_{\partial _{i-n+4+(-1)^i(1-4(1-\lfloor \frac{i  \ mod \ 4}{2} \rfloor))}}$ forms an $n$-path with $M_{\partial _i}$.

\end{enumerate}
\end{lemma}

\begin{proof}

This may be checked by direct inspection.

\end{proof}


Now we focus our attention on rows $n \equiv 1 \ mod \ 4$. From
Lemma \ref{lema1}, it is clear that $n$-paths consists of groups
of coboundaries in the same coset modulo 4.

\begin{lemma} \label{lema2}
Given $ 1 \leq i \neq j \leq 4t$, $i \equiv j \ mod \ 4$, there
exists one and only one row $n$, $5 \leq n \leq 2t+2$, $n \equiv 1
\ mod \ 4$, such that $M_{\partial _i}$ and $M_{\partial _j}$ form an
$n$-path.
\end{lemma}

\begin{corollary}
Along the $\frac{t-1}{2}$ rows $n \equiv 1 \ mod \ 4$, $5 \leq n \leq
2t+2$, any $k$ coboundaries $\partial _{i_1}, \ldots ,\partial
_{i_k}$ in the same coset modulo 4 give rise to a total amount of
$\frac{k(t-k)}{2}$ paths.
\end{corollary}

\begin{proof}

Along the $\frac{t-1}{2}$ rows $n \equiv 1 \ mod \ 4$, $5 \leq n \leq
2t+2$, any $k$ coboundaries $\partial _{i_1}, \ldots ,\partial
_{i_k}$ in the same coset modulo 4 might give rise to $k
\frac{t-1}{2}$ paths. Actually, this is not the case, since we
know from Lemma \ref{lema2} that every pair of such coboundaries
forms a path at one and only one of these rows $n$. Thus the total
amount of paths has to be reduced in the number of pairs in which
the $k$ coboundaries may be grouped. This gives $k \frac{t-1}{2}-
\frac{k(k-1)}{2}=k \frac{t-k}{2}$, as claimed.

\end{proof}

\begin{example}
For $t=5$, the set of $k=2$ coboundaries $\{ \partial_{14},
\partial_{18} \}$ defines exactly $2 \frac{5-2}{2}=3$ paths  in rows
congruent to $1$, namely, $1$ path at row $5$, $\{ (\partial_{18}
,  \partial_{14} ) \}$  and $2$ paths at row $9$, $\{ (
\partial_{14} ) , (
\partial_{18}) \}$.

\end{example}

Table \ref{tabla1} in Appendix \cite{AGG}, shows the total amount
$\frac{k(t-k)}{2}$ of paths produced by $k$ coboundaries in the
same coset modulo 4, for odd values of $t$.

%

\begin{proposition}\label{proposicion1} This table has many valuable combinatorial properties:

\begin{enumerate}

\item The table is symmetric.

\item The numbers in the central columns are triangular numbers, of the type $\frac{n(n+1)}{2}$.

\item Subtracting from a number in the central columns any of the numbers of the same row, gives as result a triangular number as well.

\item Reciprocally, subtracting from a number in the central columns any triangular number gives as result a number of the same row.

\end{enumerate}

\end{proposition}

\begin{proof}

\begin{enumerate}
\item The table is symmetric, since $k$ coboundaries give rise to $\frac{k(t-k)}{2}$ paths, exactly the same amount of paths produced by $t-k$
coboundaries, $\frac{(t-k)k}{2}$.
\item The numbers in the central columns are triangular numbers. Actually, $\frac{t-1}{2}$ coboundaries give rise to
$\frac{n (n+1)}{2}$ paths, where  $t=2n+1$.
\item Subtracting from a number in the central columns any of the numbers of the same row, gives as result a triangular number as well. Indeed,
subtracting $\frac{k(t-k)}{2}$ paths from $\frac{t^2-1}{8}$ gives
$\frac{t^2-1-4kt+4k^2}{8}=\frac{(t-2k)^2-1}{8}=\frac{\frac{t-2k-1}{2}\cdot
\frac{t-2k+1}{2}}{2}.$
\item The argument above fits here as well.
\end{enumerate}
\end{proof}

Attending to the condition (\ref{filahadamard1}), in order to get
a $\Z _t \times \Z_2^2$-cocyclic Hadamard matrix, a necessary (but
not sufficient!) condition is to select $k_i$ coboundaries in the
coset $i \ mod \ 4$, such that there is a total amount of
$t\frac{t-1}{2}$ paths along the $\frac{t-1}{2}$ rows $n \equiv 1
\ mod \ 4$, $5 \leq n \leq 2t+2$. This motivates the following
definition.

\begin{definition}\label{definicion1}
A {\em distribution} is a tuple $(\frac{k_0(t-k_0)}{2},\frac{k_1(t-k_1)}{2},\frac{k_2(t-k_2)}{2},\frac{k_3(t-k_3)}{2})$, $0 \leq k_j
\leq k_i \leq \frac{t-1}{2}$ for $j \geq i$, such that
\begin{equation}\label{ecuacionfila1} \sum _{i=0}^3
\frac{k_i(t-k_i)}{2}=\frac{t(t-1)}{2}.\end{equation}
\end{definition}

\begin{proposition} \label{proposicion2}
For any odd $t$, there always exists at least one distribution
$(\frac{k_0(t-k_0)}{2},\frac{k_1(t-k_1)}{2},\frac{k_2(t-k_2)}{2},\frac{k_3(t-k_3)}{2})$.
Furthermore, there are as many different distributions as
decompositions of $\frac{t-1}{2}$ as the summation of four
triangular numbers.
\end{proposition}

\begin{proof}

The maximum possible number of paths is $4
\frac{t^2-1}{8}=\frac{t^2-1}{2}$,  when $k_i= \frac{t-1}{2},
\forall i$, so that the relation (\ref{ecuacionfila1}) fails to
hold by a difference $m=\frac{t^2-1}{2}-\frac{t(t-1)}{2}$.

In 1796 Gauss proved that any positive integer can be decomposed
as the summation of three (not necessarily different) triangular
numbers, some of which may be eventually zero. Consequently, there
exist three triangular numbers $0 \leq t_1,t_2,t_3 \leq
\frac{t^2-1}{8}$ such that $m=t_1+t_2+t_3$.

Thus $\frac{t^2-1}{2}=4
\frac{t^2-1}{8}-m=\frac{t^2-1}{8}+(\frac{t^2-1}{8}-t_1)+(\frac{t^2-1}{8}-t_2)+(\frac{t^2-1}{8}-t_3)$.
Taking into account Proposition \ref{proposicion1}, there exist
integers $0 \leq k_3 \leq k_2\leq k_1 \leq \frac{t-1}{2}$ such
that  $(\frac{t^2-1}{8}-t_i)=\frac{k_i(t-k_i)}{2}$, and therefore
$(\frac{t^2-1}{8},\frac{k_1(t-k_1)}{2},\frac{k_2(t-k_2)}{2},\frac{k_3(t-k_3)}{2})$
is a distribution, in the sense of Definition \ref{definicion1}.

The second part is a straightforward consequence.

\end{proof}

This proposition provides a method for finding the set of
distributions for a given $t$, in terms of decompositions of
$\frac{t-1}{2}$ as the summation of four triangular numbers $0
\leq t_0 \leq t_1\leq t_2 \leq t_3 \leq \frac{t^2-1}{8}$.

\begin{proposition} \label{generadistribucion}
Let $k$ be a positive integer. Then:
\begin{enumerate}

\item $k$ is a triangular number if and only if $\frac{-1+\sqrt{1+8k}}{2}$ is an integer.

\item The greatest triangular number less or equal to $k$ is $t_n$, for $n=\lfloor \frac{-1+\sqrt{1+8k}}{2}\rfloor$.

\item If $k$ is decomposed as the summation of $m$ triangular numbers $t_{i_j}$, $1\leq j \leq m$, then $\D \max _j \{t_{i_j}\} \geq t_n $, for $n=\lceil \frac{-1+\sqrt{1+8\frac{k}{m}}}{2} \rceil$.

\end{enumerate}
\end{proposition}

\begin{proof}

It suffices to notice that $k$ is a triangular number if and only
if there exists an integer $n$ such that $t_n=n \frac{n+1}{2}=k$.
Equivalently, if and only if the equation
$\frac{n^2}{2}+\frac{n}{2}-k$ has a positive integer solution
(which, a fortiori, is $\frac{-1+\sqrt{1+8k}}{2}$).

\end{proof}

Proposition \ref{generadistribucion} leads straightforwardly to an
algorithm for constructing the full set of distributions for a
given $t$ (see Algorithm 1 in Appendix \cite{AGG}).

%
%

Table \ref{tabla2} in Appendix \cite{AGG}, shows the complete set of
distributions obtained from Algorithm 1, for $3\leq t \leq 25$.

Notice that the knowledge of the full set of distributions implies
the knowledge about the number of coboundaries which have to be
used in order to get a $\Z _t \times \Z_2^2$-cocyclic Hadamard
matrix, since each summand $\frac{k_i(t-k_i)}{2}$ is in one to one
correspondence to the values $k_i$ and $t-k_i$ (see Table
\ref{tabla1}). In spite of this fact, we may bound the number of
coboundaries to be combined a bit further.

\begin{lemma}\label{distri}
For any decomposition of $n$ into $k$ summands, say $n=n_1+ \dots
+ n_k$ We have $ \sum _{i=1}^k n_i^2 \geq \frac{1}{k} \left( \sum
_{i=1}^k n_i\right) ^2$.
 In fact, one could check that
\begin{equation} \sum _{i=1}^k n_i^2-\frac{1}{k} \left( \sum _{i=1}^k n_i\right) ^2 = \frac{1}{k} \sum _{1 \leq j <i \leq k} (n_i-n_j)^2\geq 0. \end{equation}
\end{lemma}

\begin{proposition} \label{proposicion4}
Let
$(\frac{k_0(t-k_0)}{2},\frac{k_1(t-k_1)}{2},\frac{k_2(t-k_2)}{2},\frac{k_3(t-k_3)}{2})$
be a distribution. Call $n=k_0+k_1+k_2+k_3$. Then

\begin{enumerate}

\item
\begin{equation}\label{cotacob1} \lceil \frac{t-\sqrt{4t-3}}{2}\rceil \leq k_3 \leq \lfloor \frac{t+\sqrt{4t-3}}{2}\rfloor . \end{equation}

\item
\begin{equation}\label{cotacob2} \lceil 2(t-\sqrt{t})\rceil \leq n \leq \lfloor 2(t+\sqrt{t})\rfloor . \end{equation}

\end{enumerate}
\end{proposition}

\begin{proof}

Let $(k_0,k_1,k_2,k_3)$ generate a distribution
$(\frac{k_0(t-k_0)}{2},\frac{k_1(t-k_1)}{2},\frac{k_2(t-k_2)}{2},\frac{k_3(t-k_3)}{2})$.


On one hand,  condition (\ref{ecuacionfila1}), gives
$4k_3^2-4k_3t+t^2-4t+3 \leq 0 $ which proves (\ref{cotacob1}).

On the other hand, simplifying (\ref{ecuacionfila1}), we get $\D t
\sum _{i=0}^3k_i-t^2+t=\sum _{i=0}^3k_i^2.$ Now, by Lemma
\ref{distri}
\begin{equation}\label{repartohom} \D \sum _{i=0}^3k_i^2 \geq \sum _{i=0}^3(\frac{n}{4})^2, \end{equation} and so  $\D tn-t^2+t \geq \frac{n^2}{4}
$, obtaining  (\ref{cotacob2}).


\end{proof}

\begin{remark} \label{nota1}
Condition $(\ref{cotacob2})$ may be tightened, depending on the
coset of $n=k_0+k_1+k_2+k_3$ modulo 4, substituting the lower
bound in (\ref{repartohom}) by the most homogeneously distributed
partition of $n$ into four parts.

%
%
%
%
%
%

\end{remark}

\begin{remark}

The bounds in Proposition \ref{proposicion4} are very tight, as
has been checked experimentally. The first gap occurs for $t=71$,
and consists in just one coboundary.

\end{remark}

Once we know that a distribution is available for a given value of
$t$, the next step is looking for appropriate subsets of $n_i$
coboundaries in the cosets $i \ mod \ 4$ in ${\cal B}$ such that
the amount of $n$-paths along rows $n \equiv 1 \ mod \ 4$, $5 \leq
n \leq 2t-1$ fits that distribution.

\begin{example}
For $t=5$, the distribution $(3,3,2,2)$ corresponds to a
repartition of $(2,2,1,1)$ coboundaries in each of the cosets
modulo 4 (although we can consider $3$ instead of $2$ and $4$
instead of $1$ coboundaries in each case). Now, for each subset of
$k_i$ coboundaries, we can compute the number of $n$-paths defined
along the rows congruent to $1$ (this can be computed in only one
 coset, they are disjoint).

\end{example}

\begin{definition}\label{ingredient} An {\em ingredient} produced by a subset of $k$ coboundaries in ${\cal B}$ in the same coset modulo 4 is the column vector
whose entries are the number of $n$-paths produced by these $k$
coboundaries along rows $n \equiv 1 \ mod \ 4$, $5 \leq n \leq
2t-1$. A {\em recipe} is a collection of 4 ingredients (one for
each different coset $i$ modulo 4), arranged as a matrix of 4
columns, such that the sum of each of the rows is $t$.
\end{definition}

\begin{example}
For $t=5$, the set of $k=2$ coboundaries $\{ \partial_{14},
\partial_{18} \}$ (and also any set $\{ \partial_{14+i}, \partial_{18+i} \} \ mod \ 4t$
for any $i$) define the ingredient $[1,2]^t$.
One recipe for $t=5$, will be $\left[ \begin{array}{c} 1 \\ 2
\end{array} \right], \left[
\begin{array}{c} 2 \\ 1 \end{array} \right], \left[ \begin{array}{c} 1 \\  1 \end{array}
\right], \left[ \begin{array}{c}  1 \\ 1 \end{array} \right]$.

\end{example}

Consequently, if a subset of $\{n_1,n_2,n_3,n_4\}$ coboundaries in
${\cal B}$ defines a recipe, this subset of coboundaries satisfies
the condition (\ref{filahadamard1}), and therefore the summation
of each of the rows $n \equiv 1 \ mod \ 4$, $5 \leq n \leq 2t-1$ is
zero.

\begin{proposition}
The notion of recipe does not depend on the order of its
ingredients.
\end{proposition}

\begin{proof}

Attending to Lemma \ref{lema1}, $\partial _i$ forms an $n$-path
with $\partial _{i-n+1}$, independently on the coset $i \ mod \
4$, for $n \equiv 1 \ mod \ 4$, $5 \leq n \leq 2t-1$. In
particular, $n$-paths are constructed from those coboundaries in
${\cal B}$ in the same coset $\ mod \ 4$, which differ in
$\frac{n-1}{4}$ positions in the 5-cycles.

This way, if a subset of coboundaries of the coset $i \ mod \ 4$
produces an ingredient, the same is produced by the translation of
this subset to any other coset $\ mod \ 4$.

Eventually, this translation could produce a coboundary $\partial
_i$ which is not in ${\cal B}$. This is not a source of
difficulties, since 
such prohibited subsets of
coboundaries may be substituted by their complements in the
5-cycles. Since the substitution of any amount of paths by their
complementary in a cycle does not change the total amount of
paths, this operation preserves the ingredient.

\end{proof}

Finding a recipe is the first step in the process of constructing
a $\Z _t \times \Z_2^2$-cocyclic Hadamard matrix, since any subset
of coboundaries in ${\cal B}$ defining a recipe satisfies
condition (\ref{filahadamard1}) and conversely.

Turning our attention to rows not congruent to 1, where the number
of paths must be equal to the number of intersections in order to
fullfill the Hadamard test, the following proposition gives a
condition about when the relation (\ref{filahadamard2}) is also
satisfied (and hence a $\Z _t \times \Z_2^2$-cocyclic Hadamard
matrix has been found).

\begin{proposition}
A subset of coboundaries in ${\cal B}$ satisfies condition
(\ref{filahadamard2}) (i.e. the summation of the $n^{th}$-row is
zero, for $n\equiv 0,2,3 \ mod \ 4$, $6 \leq n \leq 2t+2$),  if
and only if the number of $n$-paths of even length  is itself
even, half of them starting and ending with coboundaries in cosets
$i_1,i_2 \ mod \ 4$, the other half starting and ending in
coboundaries in cosets $i_3,i_4 \ mod \ 4$, $i_j \neq i_k$ for $j
\neq k$.
\end{proposition}

\begin{proof}

As we commented in Section \ref{2}, we are using $R=1_t \otimes
 \left(\begin{array}{rrrr}1&1&1&1\\1&-1&1&-1 \\1&-1&-1&1\\1&1&-1&-1
\end{array}\right)$ as the matrix coming from representative cocycles.

Consequently, intersections in rows $n \equiv 2 \ mod \ 4$ can
occur in positions $(n,i)$, for $i \equiv 2,0 \ mod \ 4$.
Similarly, intersections in rows $n \equiv 3 \ mod \ 4$ can occur
in positions $(n,i)$, for $i \equiv 2,3 \ mod \ 4$. Finally,
intersections in rows $n \equiv 0 \ mod \ 4$ can occur in
positions $(n,i)$, for $i \equiv 3,0 \ mod \ 4$. Taking into
account Lemma \ref{lema1}, it follows that $n$-paths consists in
properly alternating coboundaries in:

\begin{itemize}

\item Either cosets $(1,2) \ mod \ 4$, either cosets $(3,0) \ mod \ 4$, for $n \equiv 2 \ mod \ 4$.

\item Either cosets $(2,0) \ mod \ 4$, either cosets $(1,3) \ mod \ 4$, for $n \equiv 3 \ mod \ 4$.

\item Either cosets $(2,3) \ mod \ 4$, either cosets $(1,0) \ mod \ 4$, for $n \equiv 0 \ mod \ 4$.

\end{itemize}

Hence any $n$-path of odd length produces exactly one intersection
(i.e. shares exactly one negative entry) with $R$ at the
$n^{th}$-row. On the other hand, $n$-paths of even length produces
either 2 or 0 intersections, depending on the cosets modulo 4 of
$n$ and the initial coboundary of the $n$-path. More precisely:

\begin{itemize}

\item If $n\equiv 2 \ mod \ 4$, then an $n$-path of even length will produce two intersections at the $n^{th}$-row if and only if the coset $i$ of the
initial coboundary is $i \equiv 2,0 \ mod \ 4$.

\item If $n\equiv 3 \ mod \ 4$, then an $n$-path of even length will produce two intersections at the $n^{th}$-row if and only if the coset $i$ of the
initial  coboundary is $i \equiv 2,3 \ mod \ 4$.

\item If $n\equiv 0 \ mod \ 4$, then an $n$-path of even length will produce two intersections at the $n^{th}$-row if and only if the coset $i$ of the
initial  coboundary is $i \equiv 3,0 \ mod \ 4$.

\end{itemize}

Summing up,  each $n$-path of odd length produces 1 intersection,
and each $n$-path of even length  produces either 2 or 0
intersections. Hence, the only circumstance in which the amounts
of intersections and $n$-paths both coincide is precisely when
half the  $n$-paths of even length give rise to 2 intersections,
whereas the remaining half of  $n$-paths of even length  do not
produce any intersections at all.

\end{proof}

\begin{example}
For $t=5$, the set  $\{ \{
\partial_{14},
\partial_{18} \} , \{ \partial_{3}, \partial_{11} \}, \{\partial_{8} \} , \{ \partial_{5}\} \}$, of $6=2+2+1+1$ coboundaries corresponding to the recipe
showed after Definition \ref{ingredient} defines a Hadamard
matrix, because paths of even length at rows congruent to $2,3,0$
are balanced. On the other hand, the set $\{ \{
\partial_{10},
\partial_{14} \} , \{ \partial_{3}, \partial_{11} \}, \{\partial_{8} \} , \{ \partial_{5}\} \}$ (which corresponds to the same recipe) does not define a Hadamard
matrix, because it fails the Hadamard test at row $8$, being its
$8$-paths $ ( \partial_{14}, \partial_{11}) , ( \partial_{10} ) ,
( \partial_3 ), ( \partial_{5} ) , ( \partial_{8})$, there is only
one path of even length, so the number of intersections at row $8$
can not be equal to the number of $8$-paths.

\end{example}

Now it is straightforward to design an algorithm (see algorithm 2
in Appendix \cite{AGG}), searching exhaustively for $\Z _t \times
\Z_2^2$-cocyclic Hadamard matrices for odd $t$.

%
%





Table \ref{tabla3} shows an exhaustive calculation of $\Z _t
\times \Z_2^2$-cocyclic Hadamard matrices (last column) for odd
$t$, $3 \leq t \leq 13$, in terms of distributions (second
column), number of different ingredients produced by the
coboundaries needed (third column) and recipes found (fourth
column). The fifth column shows how many
of the recipes are ``productive", in the sense that they actually
give rise to $\Z _t \times \Z_2^2$-cocyclic Hadamard matrices.

\begin{table}[h]\renewcommand{\arraystretch}{1.3}\caption{$\Z _t \times \Z_2^2$-cocyclic Hadamard matrices from
Algorithm 2}\label{tabla3}\centering
$\begin{array}{c|c|c|c|c|c} t& distribution&\# ingredients & \# recipes& \# Had. recipes&\# Had. matrices\\ \hline 3&(1,1,1,0)& (1,1,1,1)&4&4&24\\
 \hline 5&(3,3,2,2)& (2,2,1,1)&12&12&120\\
 \hline 7&(6,6,6,3)& (4,4,4,1)&28&24&336\\
 &(6,5,5,5)& (4,3,3,3)&60&36&504\\
 \hline 9& (10,10,9,7) &(10,10,7,4)&756& 108&1944\\
  \hline 9&(9,9,9,9)& (7,7,7,7)&60&24&432\\
 \hline 11&(15,14,14,12)& (26,20,20,10)&5580&120&2640\\
 \hline 13& (21,21,21,15)& (74,74,74,14)& 19320 & 144&3744\\
&(21,21,18,18)& (74,74,34,34)&29208&72&1872\\
&(20,20,20,18)& (57,57,57,34)&21612&108&2808\\
 \hline
\end{array}$
\end{table}

\section{On basis, generators and diagrams}\label{5}

In this section we introduce some elementary notations, and give
the notion of {\em diagram}, a very useful presentation of a $\Z
_t \times \Z_2^2$-cocyclic matrix. In particular, the notion of
symmetry in a diagram (see Definition \ref{defsimetria}) will lead
to a fast Hadamard test for $\Z _t \times \Z_2^2$-cocyclic
matrices, which we state in Theorem \ref{teorema1}.

Cocyclic matrices over  $\Z_t \times \Z_2^2$ can be visualized by
diagrams which represent the coboundaries which appear in the
expression of the matrix.

\begin{definition}
Given a pointwise product of coboundaries  $ M_{\partial_{d_1}}
\ldots M_{\partial_{d_k}} \cdot  R$ defining a cocyclic matrix
over $\Z_t \times \Z_2^2$,  its  {\it{diagram}} is a  $4 \times t$
matrix, such that  $\{ a_{ij} \}_{ 1 \leq i \leq 4, \ 1\leq j \leq
t}$ is $\times$, ($ a_{ij} = \times$)  if  $ 4t - 4(j-1)-3 + i \
mod \ 4t \in \{ d_1, \ldots , d_k  \}$
 and empty elsewhere.
\end{definition}

\begin{remark}
The definition of diagram has to do with the expression of the
matrix in terms of the coboundaries, so every cocyclic matrix over
$\Z_t \times \Z_2^2$ has eight different diagrams.

\end{remark}

\begin{example}
For instance, one Hadamard matrix for $t=5$ is given by the
following coboundaries $\{ \{14, 18 \}, \{3, 11\}, \{ 8 \}, \{ 5
\} \}$. A presentation of this subset of coboundaries as a $4 \times t$
matrix, such that coboundaries $\ mod \ i$ are placed at row $i-1$
is:
$$\left|
    \begin{array}{ccccc}
      {\bf {18}} & {\bf {14}} & 10 & 6 & 2 \\
       19 & 15 & {\bf {11}} & 7 & {\bf {3}} \\
       20 & 16 & 12 & {\bf {8}} & 4 \\
       17 & 13 & 9 & {\bf {5}} & 1 \\
    \end{array}
  \right|, \mbox{or in short}, \left|
    \begin{array}{ccccc}
          \times & \times & - & - & - \\
        - & -  & \times & - & \times \\
        - & - & - &  \times & -  \\
        - & - & - & \times  & -\\
    \end{array}
  \right|$$
\end{example}

Diagrams are a very useful tool. This presentation of the
coboundaries allows us to read easily  the adjacency conditions,
the number of paths, their length and the number of intersections
they produce at each row $n$.

Every row of a diagram represents one  $5$-cycle.
 Actually, the diagram is not a matrix, but a cylinder, the
first and last columns being adjacent. The paths in rows  not
congruent to $1$ can be read by alternating pair of rows in the
diagram. In the previous example, the set of coboundaries selected
defines the following $8$-paths: $ (14 , 11) , ( 3, 18 ) , ( 8 ) ,
( 5 ).$

Among the two paths of even length, $( 3, 18 ) $ produces 2
intersections (at positions $3$ and $15$), and $( 14 , 11 )$
produces no intersection (it has $-1$ at positions $14$ and $6$).
So whether they give $0$ or $2$ intersections comes from the
congruency class module $4$ of its first (or last) coboundary.
This can be easily checked in any diagram.

\begin{definition}\label{defsimetria}
A diagram associated to a pointwise product of coboundaries  $M =
M_{\partial_1} \ldots M_{\partial_k} R$ is called {\it{symmetric}}
if the $\times$ are symmetric with respect to a column. If one of
the diagrams representing a cocyclic matrix is symmetric, so are
all other diagrams (they are obtained by complementing some of the
rows, and symmetry is preserved). We will say, by extension, that
the cocyclic matrix presents symmetry (the matrix is not
symmetric, the diagram is).
\end{definition}

\begin{theorem} \label{teorema1}
If a cocyclic matrix over   $\Z_t \times \Z_2^2$ is represented by
a symmetric diagram and has exactly $t$ paths in rows congruent to
$1 \ mod \  4$, then it is Hadamard.
\end{theorem}

\begin{proof}

Taking into account Proposition \ref{corolario1}, we have to prove
that the number of paths coincide with the number of intersections
for rows $6$, $7$, $8$ and congruent. Each path of odd length
give one intersection, and the only thing to prove is that  paths
of even length  are equally distributed between those giving $0$
intersections and those giving $2$ intersections. As these paths
of even length come in pairs (they do not use any of the
coboundaries on the symmetry axis) the character of each path is
the opposite of its symmetric, so its number is balanced and the
result holds.

\end{proof}

\section{Operations}\label{6}

In this section we will study four different operations on the set
of coboundary matrices over $\Z_t \times \Z_2^2$: complements,
rotations, swappings and dilatations. In particular, these
operations extend to operations over $\Z_t \times \Z_2^2$-cocyclic
matrices, which will be proved to preserve orthogonality. These operations partition the set of $\Z_t
\times \Z_2^2$-cocyclic matrices into disjoint orbits, which can
be easily computed once one element is known. Among all the
elements of an orbit, a representative can be chosen, in a
standard way that will be made precise.

\subsection{Complements. The group $\Z_2$}

Given a set of coboundaries, $t$, and one of the four congruency
classes  $\ mod \  4$, for instance $i$-th, $1 \leq i \leq 4$, we can
consider the cocyclic matrix defined by the set of coboundaries
obtained by substituting the subset of coboundaries belonging to
the $i$-th congruency by its complement.

\begin{definition}
Let   $ \{ \{c_{2 j_2} \},\{ c_{3 j_3}\},\{c_{4 j_4}\}, \{c_{1
j_1} \} \}$ with $j_k \equiv k \ mod \   4$, be a set of
coboundaries. The {\it{ complement in the $i$-th component }}, $1
\leq i \leq 4$, denoted $ {\textsc{C}_i} (\{ \{c_{2 j_2} \},\{
c_{3 j_3}\},\{c_{4 j_4}\}, \{c_{1 j_1} \} \} )$, is the union of
the complement of  $\{ c_{i j_i}\}$ in the  $i$-th component and
the rest of the initial coboundaries.
\end{definition}

For instance , if we choose  $i=2$:
$$  {\textsc{C}_2}  ( \{ \{14, 18\}, \{3, 11 \}, \{8\}, \{5\} \}) =
 \{ \{2, 6, 10 \}, \{3, 11 \}, \{8\}, \{5\} \}
$$
%
%

\begin{lemma} The only complement to compute is the
complement with respect to the component congruent to $2$.
\end{lemma}

\begin{proof}

If we consider any other complement, there will appear one of the
coboundaries  $1$, $4t-1$ or $4t$. Substitution of the expression
of this coboundary in terms of the basis of coboundaries gives us the complement with respect to the
congruency $2$.

\end{proof}

If two different sets of coboundaries give the same cocyclic
matrix their complements define the same matrix. It suffices to
compute the image of the  different  expressions for a cocyclic
matrix and check that we get the  different expressions for its
image. So, there is no imprecision when we say the {\it complement
matrix} of a cocyclic matrix over $\Z_t \times \Z_2^2$.

\begin{theorem}
 If a set a coboundaries define a Hadamard matrix, so does its complement.
\end{theorem}

  \begin{proof}

It suffices to check that the hypotheses of Theorem \ref{teorema1}
are satisfied.

In rows congruent to $1$, the complement of a set of coboundaries
in a cycle is another set which determines the same number of
paths. For the other rows, just observe that the complement
preserves symmetry in the diagram.

\end{proof}

\begin{remark}
Note that the operation {\it complement} modifies the number of
coboundaries, substituting  the $k_0$ coboundaries congruent to $2
\ mod \   4$ by $t-k_0$ coboundaries.

\end{remark}

This operation can be observed as an action of the group $\Z_2$
over the set of $\Z_t \times \Z_2^2$-cocyclic matrices. For every
Hadamard cocyclic matrix over  $\Z_t \times \Z_2^2$, we can
consider its orbit. This orbit has exactly $2$ elements.

\subsection{Rotations. The group $\Z_t$}

Given that the Hadamard condition for a cocyclic matrix over $\Z_t
\times \Z_2^2$ depends on the conditions of position and adjacency
of the coboundaries defining such matrix, one natural idea is
rotating the positions of the coboundaries.

\begin{definition}
Let   $ \{ \{c_{2 j_2} \},\{ c_{3 j_3}\},\{c_{4 j_4}\}, \{c_{1
j_1} \} \}$ with $j_k \equiv k \ mod \   4$, be a set of
coboundaries. The {\it{$i$-rotated set}},  $0 \leq i \leq t-1$ of
this set of coboundaries, denoted by $\textsc{T}_i ( \{ \{c_{2
j_2} \},\{ c_{3 j_3}\},\{c_{4 j_4}\}, \{c_{1 j_1} \} \})$ is the
set:
$$  \{ \{c_{2 j_2 -4i} \},\{c_{3 j_3 -4i }\}, \{c_{4 j_4 -4i} \}, \{ c_{1 j_1-4i}\} \} $$
where additions are $\ mod \  4t$.
\end{definition}

For instance, for  $i=2$
$$ \textsc{T}_2 ( \{ \{14, 18\}, \{3, 11\}, \{8\}, \{5 \} \} )= \{ \{ 6, 10 \}, \{3, 15 \}, \{20\}, \{ 17  \} \}$$

%

The  $i$-rotation operation moves the marked positions  $i$ places
to the right.

As in the complement operation, $i$-rotation over each of the
eight expressions for a Hadamard matrix gives each of the eight
expressions of the same cocyclic Hadamard matrix, so we are
rigourous when speaking of the  $i$-rotated of a matrix.

This operation can be observed as an action of the group $\Z_t$
over the set of $\Z_t \times \Z_2^2$-cocyclic matrices. The
element  $i$,
 $0 \leq i \leq t-1$ acts on every cocyclic matrix by substracting  $4i \ mod \   4t$
 to each of the coboundaries defining the matrix.

\begin{theorem}
The $i$-rotated set of a set defining a cocyclic Hadamard matrix
over $\Z_t \times \Z_2^2$ defines a Hadamard matrix too.
\end{theorem}

\begin{proof}


Rotations preserve the relative positions of the coboundaries in
the diagram and, thus, the number, length and intersections
defined by the paths.

\end{proof}

\begin{remark}
There is no need to consider symmetry condition for this proof.
For every Hadamard cocyclic matrix over  $\Z_t \times \Z_2^2$, we
can consider its orbit. This orbit has exactly $t$ elements.

\end{remark}

\begin{remark}
Notice that the complement and rotation operations commute.

\end{remark}

\subsection{Swappings. The group $S_4$}

The word  {\it{swapping}} explains clearly the effect of this
operation on a set of coboundaries. The marked positions are
permuted among the rows in the diagram.

Given the set of coboundaries $\{c_{i j_i}\}$ for $1 \leq i \leq
4$, we will denote $\{c_{i j_i}\}+ k$ the set of coboundaries
obtained by adding $k$ to each of their indexes.


\begin{definition}
For the set of coboundaries  $  \{ \{c_{2 j_2}\},\{c_{3
j_3}\},\{c_{4 j_4}\}, \{c_{1 j_1}\} \} $ we define the
{\it{swapping operations}}:
\begin{itemize}
\item $s_{12} ( \{ \{c_{2 j_2}\},\{c_{3 j_3}\},\{c_{4 j_4}\}, \{c_{1 j_1}\} \} ) = \{ \{c_{1 j_1}\} + 1,\{c_{3 j_3}\},\{c_{4 j_4}\}, \{c_{2 j_2}\} -1 \}. $
\item $s_{13} ( \{ \{c_{2 j_2}\},\{c_{3 j_3}\},\{c_{4 j_4}\}, \{c_{1 j_1}\} \} ) = \{ \{c_{2 j_2}\},\{c_{1 j_1}\}+2,\{c_{4 j_4}\}, \{c_{3 j_3}\} - 2 \}. $
\item $s_{14} ( \{ \{c_{2 j_2}\},\{c_{3 j_3}\},\{c_{4 j_4}\}, \{c_{1 j_1}\} \} ) = \{ \{c_{2 j_2}\},\{c_{3 j_3}\},\{c_{1 j_1}\}+3, \{c_{4 j_4}\} -3 \}. $
\end{itemize}
\end{definition}
Although we have defined three operations, any composition of them
can be considered, obtaining one of the possible $24$ permutations
on the rows of the diagram.

For instance, for  $t=5$
$$ s_{13} ( \{ \{14, 18\}, \{3, 11\}, \{8\}, \{5 \} \} )= \{ \{ 14, 18 \}, \{ 7 \}, \{8 \}, \{ 1, 9  \} \}$$

%

\begin{theorem}
If $  \{ \{c_{2 j_2}\},\{c_{3 j_3}\},\{c_{4 j_4}\}, \{c_{1 j_1}\}
\} $ is a set of coboundaries defining a Hadamard cocyclic matrix
over $\Z_t \times \Z_2^2$, then  $s_{ij} (  \{ \{c_{2
j_2}\},\{c_{3 j_3}\},\{c_{4 j_4}\}, \{c_{1 j_1}\} \} )$, $1 \leq i
< j \leq 4$ is also a Hadamard cocyclic matrix over $\Z_t \times
\Z_2^2$.
\end{theorem}

\begin{proof}

Once again, it suffices to check whether the hypotheses of Theorem
\ref{teorema1} are satisfied. On one hand, any swapping preserves
the number of paths in rows congruent to $1$. On the other hand,
any permutation of the rows in the diagram preserves symmetry.

\end{proof}

\begin{remark}
This result can be proved without using the symmetry condition,
one only need to assume that the  paths of even length giving $2$
and $0$ intersections, are equally distributed among the two
subsets of components defined in each row, attending to the
congruency
  $\ mod \   4$ of the coboundaries.

\end{remark}

This gives us an action of $S_4$ on diagrams, which commute with
complement and rotations. Depending on the diagram (i.e. on the
size of the permutation group of its rows), the orbit can have
less than $24$ elements (actually, 1,4,6,12 or 24).

\subsection{Dilatations. The group $\Z_t^*$}

As we have seen, the  orthogonality of a cocyclic matrix over
$\Z_t \times \Z_2^2$ is determined by the position and adjacency
of the coboundaries which define it. This is the reason because of
rotating the rows of a diagram preserves the orthogonality of
the corresponding $\Z_t \times \Z_2^2$-cocyclic matrix. Now we are
interested in other kind of geometric transformations:
homothecies.

\begin{definition}
Given a set of coboundaries $S$ determining a diagram with col\-umns
$ ({\cal{C}}_{t},\cdots ,{\cal{C}}_1)$, the {\it{$r$-th
dilatation}} applied to $S$, with $r \in \Z_t^*$, denoted
by $V_r (S)$, is the set corresponding to the image of the diagram
under the homothecy with ratio $r$ and center the column placed at
the right hand-side of the diagram.

As the column located at the right-hand side corresponds to the
coboundaries $2,3,4,1$, a formula for $V_r$ is given by
$V_r(\partial _k)=\partial _h$, for \begin{equation}
\label{dilatacion} h=4 \left( \left( \frac{k - (k\ mod \   4)}{4}
\cdot r\right) \ mod \ t \right) + (k \ mod \ 4). \end{equation}
\end{definition}

For instance, for  $r=2$
$$ \textsc{V}_2 ( \{ \{14, 18\}, \{3, 11\}, \{8\}, \{5 \} \} )= \{ \{ 6, 14 \}, \{3, 19\}, \{12 \}, \{ 9\} \}$$

%

\begin{remark}
Notice that a dilatation $V_r$ defines a bijection over the set of
$\Z_t \times \Z_2^2$-coboundaries if and only if $r \in \Z_t^*$.

\end{remark}

\begin{theorem}
Dilatations $V_r$ with $r \in \Z_t^*$ preserve the orthogonality of a $\Z_t \times \Z_2^2$-cocyclic matrix.
\end{theorem}

\begin{proof}

It suffices to check whether the hypotheses of Theorem
\ref{teorema1} are satisfied. In rows congruent to $1$,
dilatations define some permutation on the number of paths defined
(a permutation of the entries in the ingredient), so the total
number of paths remains equal to $t$.

Since the image diagram is symmetric as well, the result holds.

\end{proof}

\begin{remark} The size of the orbit of a symmetric $\Z_t \times \Z_2^2$-cocyclic  matrix under the action of dilatations is, at most,  ${\phi(t)}$,
$\phi$ being the Euler's totient function.

\end{remark}

\begin{remark}
Depending on the diagrams, the image under dilatations can
coincide with the image under some of the previously defined
operations; for instance the action on the previous example of the
dilatation $V_2$ coincides with the action of the composition
$s_{23}  T_4$.

\end{remark}

\subsection{Orbits and representatives}

\begin{definition}
The {\it{total orbit}} of a cocyclic Hadamard matrix over $\Z_t
\times \Z_2^2$ is the union of all orbits under the action of
complement, rotations, swappings and dilatations.
\end{definition}

\begin{definition}
A {\it representative} of a total orbit is a set of coboundaries
associated to a diagram, with minimal number of coboundaries,
symmetric with respect to the central column, and with an
increasing number of coboundaries on each row of the diagram.
\end{definition}

\begin{example}
One Hadamard matrix for $t=5$ is given by the following
coboundaries $\{ \{14, 18 \}, \{3, 11\}, \{ 8 \}, \{ 5 \} \}$, and
the representative for its total orbit will be $\{ \{10 \}, \{
11\}, \{ 8, 16 \}, \{ 1, 17 \} \}$.

\end{example}
%

Table \ref{tabla4} in Appendix \cite{AGG}, gives us the description of all
Hadamard matrices obtained in  Table \ref{tabla3}, split in
orbits, which can be generated by its representatives under the
action of the four operations previously defined. The second
column shows the size of the orbit under
complement/rotation/swapping/dilatations, when providing new
matrices.

\section{The symmetric case}\label{7}

We have computed the full set of $\mathbf{Z} _t \times
\mathbf{Z}_2^2$-cocyclic Hadamard matrices for $t \leq 13$.

There are two special repartitions of coboundaries attending to
their congruency,  for every Hadamard matrix. One is the
repartition given by the expression of the matrix with respect to
the basis. And the other is the minimal repartition (this last one
can correspond to the matrix, or to its complement). When
performing calculations we can restrict to the minimal
repartition, in order to simplify the search, avoiding 15 of the
16 possible repartitions for each distribution.

The  rotation, swapping and dilatations preserve the number of
coboundaries involved. On the other hand, the image of the
Hadamard matrix does not depend on the chosen representation (from
the possible eight). Thus we can obtain the whole orbit of the
matrix under all these operations by only knowing the expression
of the matrix with the minimal repartition, so we can restrict the
searching to this easier case.

Once this is done, we can compute the complement of each of the
obtained matrices, to complete the orbit under all operations.

\subsection{The symmetric-diagram case ($t \in [15,23]$)}

The diagrams of all the matrices obtained for $t \leq 13$ were
symmetric (symmetry is a sufficient condition to pass the Hadamard
test in rows congruent to $2,3, 0$ module $4$, so, assuming
symmetry, the Hadamard test reduces to check if the number of
paths in rows congruent to $1$ is exactly $t$).

By following this idea, we have computed the set of cocyclic
Hadamard matrices over $\Z_t \times \Z_2^2$, having a symmetric
diagram for $15 \leq t \leq 23$, which are listed in Tables
\ref{sime} and \ref{sime2} in Appendix \cite{AGG}, in terms of their
representatives and the size of their total orbit under
complement/rotation/swapping/dilata\-tions, when providing new
matrices.

\subsection{Williamson type matrices}

The first row of Table \ref{will} shows the number of Williamson
type  matrices obtained in \cite{BH95}, together with the total
number of cocyclic Hadamard matrices over $\Z _t \times \Z_2^2$
which we have computed.

\begin{table}[h]\renewcommand{\arraystretch}{1.3}\caption{Williamson type / $\Z _t \times \Z_2^2$-cocyclic Hadamard matrices (only for symmetric diagrams if  $t\geq 15$)}\label{will}\centering
$\begin{array}{c|c|c|c|c|c|c|c|c|c|c|c}
t & 3 & 5 & 7 & 9 & 11 & 13 & 15 & 17 & 19 & 21 & 23  \\
\hline \sharp Will. & 8 & 24 & 120 & 264 & 240 & 648 & 576 & - & - & - & -  \\
\sharp H.& 24 & 120 & 840 & 2376 & 2640 & 8424 & 8640 & 13056 & 34200 & 31248 & 12144\\
  \hline
\end{array}$
\end{table}

It can be observed that  the set of Williamson type Hadamard
matrices seems to be in proportion $\frac{1}{t}$ with respect to the
set of $\Z _t \times \Z_2^2$-cocyclic Hadamard matrices consisting of symmetric diagrams. We now prove the validity of this fact.

By Lemma 3.4 in \cite{BH95}, any Williamson type matrix is
Hadamard equivalent to a cocyclic matrix over $\Z _t \times
\Z_2^2$, which is $t\times t$ back-circulant by  blocks $4 \times
4$,
\begin{equation}\label{wilcoc} \left( \begin{array}{cccc} W_1 &
\ldots &  \ldots& W_t\\ W_2&_{\cdot}\cdot^{\cdot}&W_t&W_1\\ _{\cdot}\cdot^{\cdot}&
_{\cdot}\cdot^{\cdot}&_{\cdot}\cdot^{\cdot}&_{\cdot}\cdot^{\cdot}\\ W_t&W_1& \ldots&W_{t-1}
\end{array}\right), \qquad W_i=\left(
\begin{array}{rrrr} n_i&x_i&y_i&z_i\\ x_i&-n_i& z_i&-y_i\\
y_i&-z_i&-n_i&x_i\\ z_i&y_i&-x_i&-n_i \end{array}\right),\end{equation} with  $W_{i+1}=W_{t-i+1}$ for
$1 \leq i \leq t-1$.

This matrix  (\ref{wilcoc}) is a pointwise product of the
representative matrix $R$ and a certain set of generalized
coboundaries $\{M_{\partial _{d_1}}, \ldots ,
M_{\partial_{d_k}}\}$. The additional condition
$W_{i+1}=W_{t-i+1}$ for $1 \leq i \leq t-1$, leaving  $W_1$ alone
means that the symmetry column in the diagram representing the
matrix (\ref{wilcoc}) is,  precisely, the last one to the right,
which is associated to the coboundaries $\partial _2, \partial _3,
\partial _4, \partial _1$. This fact leads us to conclude that:

\begin{proposition}
Let ${\cal H}$ the set of cocyclic Hadamard matrices over $\Z_t
\times \Z_2^2$, having a symmetric diagram. Then, the set of
Williamson type matrices of type (\ref{wilcoc}) is a subset of
${\cal H}$ of size $\frac{|{\cal H}|}{t}$. Moreover, for each
element $H \in {\cal H}$, one and only one of the rotated matrices
$T_iH$, $1 \leq i \leq t$, is a Williamson type matrix.
\end{proposition}

\begin{proof}

The condition needed tells us that the diagram has to be symmetric
with respect to the last column, which is only possible for one of
the rotated.

\end{proof}

So, the predicted number of Williamson type matrices for $t \in
[17,23]$ will be $13056/17=768$, $34200/19=1800$, $31248/21=1488$
and $12144/23=528$.

\subsection{More Hadamard matrices ($t \in [25,63]$)}

Moreover, the result about Williamson type matrices allows us to
go beyond $t=23$. Actually, we have identified every Williamson
type  matrix exhibited in the Koukouvinos website \cite{Kou} (which corresponds to a
certain decomposition of $4t$ as a sum of squares) as a pointwise
product of coboundaries and computed the whole orbit of matrices
for $t \in [25, 39]$. We show the results obtained in Table
\ref{kou} in Appendix \cite{AGG}, and give the size of its orbit under the
action of complement/rotation/swapping/di\-latations, when providing
new matrices.

The Koukouvinos website \cite{Kou} only gives an example of each
equivalence class of Williamson type Hadamard matrices for small orders, and not the total
number of these matrices. However we have checked that the numbers
of Williamson type Hadamard matrices we predict by using our computation of
symmetric-diagram cocyclic Hadamard matrices and dividing by $t$, for $t\in [17, 23]$,
is the same as if we take the examples in \cite{Kou} in each order
and compute the sum of their orbit sizes under our operations.

We conjecture that one can obtain the total number of
Williamson type Hada\-mard matrices for $t\in [25,63]$ from \cite{Kou} by just computing the sum of the orbit
sizes for each order. That is to say, we conjecture that any Williamson type matrix $H'$ Hadamard equivalent to a given Williamson type matrix $H$ is in the orbit of $H$ under our operations.

Finally we show the total number of Hadamard matrices that we have
computed so far

\begin{table}[h]\renewcommand{\arraystretch}{1.3}\caption{$\Z _t \times \Z_2^2$-cocyclic Hadamard matrices computed}\label{total}\centering
$\begin{array}{c|c|c|c|c|c|c|c|c|c|c|c}
t & 3 & 5 & 7 & 9 & 11 & 13 & 15 & 17 & 19 & 21 & 23  \\
\sharp H.& 24 & 120 & 840 & 2376 & 2640 & 8424 & 8640 & 13056 & 34200 & 31248 & 12144\\
  \hline
t  & 25 & 27 & 29 & 31 & 33 & 35 & 37 & 39 & 41 & 43 & 45  \\
\sharp H.& 75000 & 64152 & 19488 & 33480 & 79200 & - & 53280 & 7488 & 4920 & 14448 & 12960\\
  \hline
  t & 47 & 49 & 51 & 53 & 55 & 57 & 59 & 61 & 63 &  &  \\
\sharp H.& - & 24696 & 58752  & - & 26400 & 24624 & - & 21960 & 9072 &  & \\
  \hline
\end{array}$
\end{table}

\subsection{Further questions}

After an observation process, some questions come into our minds:

\begin{itemize}

\item All $\Z _t \times \Z_2^2$-cocyclic Hadamard matrices obtained so far have a symmetric diagram up to $t=13$, and we have taken advantage of it
when computing them for $t \in [15,23]$ by restricting ourselves
to the symmetric case. Can this symmetry assumption be proved? Does a $\Z _t \times \Z_2^2$-cocyclic Hadamard matrix exist coming from  a non symmetric diagram?

\item In case that the answer to the second question above is affirmative, then it would be interesting to look for $\Z _t \times \Z_2^2$-cocyclic Hadamard matrices in those orders for which no Williamson type Hadamard matrix is known to exist, such as $t=35$ \cite{Kou}.

\item Now that we have proved that our operations preserve orthogonality, arises the question: do these operations preserve
Hadamard equivalence classes? Rotations, dilatations and some of the
swappings do, because they can be expressed in terms of the bundle
operation defined by Horadam in \cite{Hor07} (this fact will be
detailed elsewhere). In addition, for values of $t$ below 23, the
number of orbits coincide with the number of non equivalent
Williamson type matrices computed in \cite{Kou}, so  the result is
probably true in general.

\end{itemize}


\begin{acknowledgements}

The authors want to express their gratitude to the anonymous referees for their comments and suggestions, which have permitted to improve the readability
and understandability of the paper.

All authors are partially supported by the research project FQM-016 from JJAA (Spain).

\end{acknowledgements}


\newpage


\appendix{Tables and algorithms}

\begin{algorithm}\label{algoritmo1} Algorithm 1: Constructing the set of distributions for $t$.\\

{\small
\noindent {\sc Input:} $t$\\
\\
$k=\{ \; \}$\\
for $i_1$ from $\lceil \frac{-1+\sqrt{t}}{2}\rceil $ to $\lfloor \frac{-1+\sqrt{4t-3}}{2} \rfloor$ do\\
\hspace*{.25cm}for $i_2$ from $\lceil \frac{-1+\sqrt{1+4\frac{t-1-i_1(i_1+1)}{3}}}{2}\rceil $ to\\
\hspace*{.25cm}$\min (i_1,\lfloor \frac{-1+\sqrt{1+4(t-1-i_1(i_1+1))}}{2} \rfloor )$ do\\
\hspace*{.5cm}for $i_3$ from $\lceil \frac{-1+\sqrt{1+2(t-1-i_1(i_1+1)-i_2(i_2+1))}}{2}\rceil $\\
 \hspace*{.5cm} to $\min (i_2,\lfloor \frac{-1+\sqrt{1+4(t-1-i_1(i_1+1)-i_2(i_2+1))}}{2} \rfloor )$\\
\hspace*{.5cm}  do\\
\hspace*{.75cm}  $x=\frac{-1+\sqrt{1+4(t-1-i_1(i_1+1)-i_2(i_2+1)-i_3(i_3+1))}}{2}$\\
\hspace*{.75cm} if $x$ is an integer, then $k= k \cup \{ (x,i_3,i_2,i_1) \}$ fi\\
\hspace*{.5cm} od\\
\hspace*{.25cm} od\\
od\\
$l=\{\;\}$\\
for $i$ from 1 to $length(k)$ do\\
\hspace*{.25cm} $l=l \cup \{ (\frac{t^2-1}{8}-k_{i,1},\frac{t^2-1}{8}-k_{i,2},\frac{t^2-1}{8}-k_{i,3},\frac{t^2-1}{8}-k_{i,4}) \}$\\
od\\
\\
{\sc Output:} $l$\\
}

\end{algorithm}

\newpage

\begin{algorithm} \label{algoritmo2} Algorithm 2: Exhaustive search for $\Z _t \times \Z_2^2$-Hadamard matrices.\\
 {\small
{\sc Input:} $t$.\\
\\
Calculate the valid distributions for $t$.\\
Calculate all the ingredients associated to every distribution.\\
Construct the set of recipes corresponding to each distribution.\\
Determine those subsets of coboundaries defining a recipe.\\
Check whether these subsets satisfy the balanced distribution of even $n$-paths for $n  \equiv 2,3,0 \ mod \ 4$, $6 \leq n \leq 2t+2$.\\
{\sc Output:} the full set of $\Z _t \times \Z_2^2$-Hadamard
matrices. }

\end{algorithm}

\renewcommand{\thetable}{\Alph{table}}

Table \ref{tabla1} shows the total amount $\frac{k(t-k)}{2}$ of
paths produced by $k$ coboundaries in the same coset modulo 4, for
odd values of $t$.

\begin{center}
\begin{table}[h]\caption{Paths produced from $k$ coboundaries in rows $n \equiv 1 \ mod \ 4$.}\label{tabla1}
$\begin{array}{c|cccccc|ccccccc} \hline t \backslash k & 1&&\dots
&&&\frac{t-1}{2}&\frac{t+1}{2}&&\dots &&&1 \\ \hline
3 && &&&&1&1&&&&&& \\\hline 5 & &&&&2&3&3&2&&&&\\ \hline 7 && &&3&5&6&6&5&3&&&\\ \hline 9 &&&4&7&9&10&10&9&7&4&&\\
\hline 11 & & 5&9&12&14&15&15&14&12&9&5&\\  \hline 13 & 6 & 11
&15&18&20&21&21&20&18&15&11&6 \\  \hline
\vdots&&&&&\vdots&\vdots&&&&&\vdots\\ \hline \end{array}$
\end{table}
\end{center}

\newpage

Table \ref{tabla2} shows the complete set of distributions
obtained from Algorithm  1, for $3\leq t \leq 25$.

\begin{table}[h]\caption{Distributions in terms of decompositions of $\frac{t-1}{2}=t_0+t_1+t_2+t_3$.}\label{tabla2}\centering
$\begin{array}{c|c|c|c} t&\frac{t(t-1)}{2}& distributions&
t_0+t_1+t_2+t_3=\frac{t-1}{2}\\ \hline 3&3&(1,1,1,0)& 0+0+0+1=1\\
\hline 5&10&(3,3,2,2)&0+0+1+1=2\\ \hline 7&21& \begin{array}{c}
(6,6,6,3)\\ (6,5,5,5)\end{array}& \begin{array}{c} 0+0+0+3=3\\
0+1+1+1=3\end{array}\\ \hline 9&36& \begin{array}{c} (10,10,9,7)\\
(9,9,9,9)\end{array}& \begin{array}{c} 0+0+1+3=4\\
1+1+1+1=4\end{array}\\ \hline 11&55&  (15,14,14,12)& 0+1+1+3=5\\
\hline 13&78& \begin{array}{c} (21,21,21,15)\\ (21,21,18,18)\\
(20,20,20,18)\end{array}& \begin{array}{c} 0+0+0+6=6\\ 0+0+3+3=6\\
1+1+1+3=6\end{array}\\ \hline 15&105& \begin{array}{c}
(28,28,27,22)\\ (28,27,25,25)\end{array}& \begin{array}{c}
0+0+1+6=7\\ 0+1+3+3=7\end{array}\\ \hline 17&136& \begin{array}{c}
(36,35,35,30)\\ (35,35,33,33)\end{array}& \begin{array}{c}
0+1+1+6=8\\ 1+1+3+3=8\end{array}\\ \hline 19&171& \begin{array}{c}
(45,45,42,39)\\ (45,42,42,42)\\ (44,44,44,39)\end{array}&
\begin{array}{c} 0+0+3+6=9\\ 0+3+3+3=9\\ 1+1+1+6=9 \end{array}\\
\hline 21&210& \begin{array}{c} (55,55,55,45)\\ (55,54,52,49)\\
(54,52,52,52)\end{array}& \begin{array}{c} 0+0+0+10=10\\
0+1+3+6=10\\ 1+3+3+3=10\end{array}\\ \hline 23&253&
\begin{array}{c} (66,66,65,56)\\ (65,65,63,60)\end{array}&
\begin{array}{c} 0+0+1+10=11\\ 1+1+3+6=11\end{array}\\ \hline
25&300& \begin{array}{c} (78,78,72,72)\\ (78,77,77,68)\\ (78,75,75,73)\\ (75,75,75,75)\end{array}& \begin{array}{c} 0+0+6+6=12\\ 0+1+1+10=12\\ 0+3+3+6=12\\ 3+3+3+3=12\end{array}\\
 \hline
\end{array}$
\end{table}

\begin{table}[h]\renewcommand{\arraystretch}{1.3}\caption{Orbits of $\Z _t \times \Z_2^2$-cocyclic Hadamard matrices}\label{tabla4}\centering
$\begin{array}{c|c|c} t& orbit & representative \\ \hline
3& 2 \times 3 \times 4 \times 1 = 24 & \{ \{ \}, \{ 7\}, \{ 8\}, \{ 5\} \}\\
 \hline 5& 2 \times 5 \times 12 \times 1 = 120  & \{ \{ 10\}, \{ 11\}, \{ 8, 16\}, \{ 1, 17\} \}\\
 \hline 7& 2 \times 7 \times 24 \times 1 = 336 & \{ \{ 14 \}, \{ 11, 15, 19 \} ,  \{ 8, 16, 24 \} , \{ 1, 13, 25 \} \} \\
\cline{2-3}   & 2 \times 7 \times  12 \times 3 = 504 & \{ \{ 10, 18 \}, \{ 11, 19 \} ,\{ 4,  28 \} , \{ 1, 13, 25 \} \}\\
 \hline 9& 2 \times 9  \times   12 \times 3 = 648 &\{ \{ 14, 22 \}, \{ 15, 19, 23 \} , \{ 4, 16, 24, 36 \} , \{ 1, 13, 21, 33 \} \}\\
\cline{2-3}  & 2 \times 9 \times 24 \times 3 = 1296 & \{ \{ 14, 22 \}, \{ 3, 19, 35 \} , \{ 12, 16, 24, 28 \} , \{ 1, 9, 25, 33 \} \}\\
\cline{2-3}  & 2 \times 9 \times 24 \times 1 = 432 & \{ \{ 14, 18,  22 \}, \{ 11, 19, 27 \} ,  \{ 8, 20, 32 \} , \{ 1, 17, 33 \} \}\\
 \hline 11& 2 \times 11 \times & \{ \{ 18, 22, 26 \}, \{ 7, 15, 31, 39 \} , \\ & 24 \times 5 = 2640& \{ 4, 16, 32, 44 \} , \{  9, 13, 21, 29, 33 \} \}\\
 \hline 13& 2 \times 13 \times & \{ \{  22, 26, 30 \}, \{ 3, 15, 23, 31, 39, 51\} , \\ & 24 \times 3 = 3744& \{ 8, 16, 20, 36, 40, 48 \} ,
\{ 1, 5, 17, 33, 45, 49 \} \}\\
\cline{2-3}  & 2 \times 13 \times & \{ \{ 18, 22, 30, 34 \}, \{ 7,
15, 39, 47 \} , \\ & 24 \times 3 = 1872 & \{ 8, 12,  24, 32, 44,
48\} ,
\{ 1, 5, 21, 29, 45, 49 \} \} \\
\cline{2-3}  & 2 \times 13 \times & \{ \{ 14, 18, 34, 38 \}, \{
15, 23, 27, 31, 39 \} , \\ & 24 \times 3 = 1872& \{ 8, 24, 28, 32,
48 \} ,
\{ 5, 17, 25, 33, 45 \} \} \\
\cline{2-3}  & 2 \times 13 \times & \{ \{ 14, 18, 34, 38 \}, \{
15, 19,  27, 35, 39 \} , \\ & 12 \times 3 = 936 & \{ 4, 12, 28,
44, 52 \} ,
\{ 1, 9, 25, 41, 49 \} \}\\
 \hline
\end{array}$
\end{table}

\begin{table}[h]\renewcommand{\arraystretch}{1.3}\caption{$\Z _t \times \Z_2^2$-cocyclic Hadamard matrices [15-19] (only for symmetric diagrams)}\label{sime}\centering
$\begin{array}{cc} t/ orbit / \# H. &  representative \\
 \hline 15 \qquad 2 \times 15 \times 24 \times 4  \qquad 2880  &  \{ \{ 18, 22, 38, 42 \}, \{ 7, 19, 27, 35, 43, 55 \} , \\  &  \{  16, 20, 28, 32, 36, 44, 48 \} , \{ 1, 9, 25, 29, 33, 49, 57 \} \} \\
\hline  15 \qquad 2 \times 15 \times 24 \times 4   \qquad 2880 &  \{ \{ 10, 14, 46, 50 \}, \{ 7, 19, 27, 35, 43, 55 \} , \\  & \{ 8, 12, 24, 32,40, 52, 56 \} , \{ 1, 5, 9, 29, 49, 53, 57 \} \} \\
\hline
 15 \qquad 2 \times 15 \times 24 \times 2  \qquad 1440 & \{ \{ 14, 26, 34, 46 \}, \{ 11, 15, 27, 35, 47, 51 \} , \\ & \{ 16,20, 28, 32, 36, 44, 48 \} , \{ 1, 5, 21, 29, 37,53, 57 \} \}\\
\hline 15 \qquad 2 \times 15 \times 12 \times 4  \qquad 1440 & \{ \{ 14, 26, 30, 34, 46 \}, \{ 15, 27, 31, 35, 47 \} , \\ & \{ 12, 16, 24, 40, 48, 52 \} , \{ 9, 13, 21, 29, 37, 45, 49 \} \}\\
\hline
\hline 17 \qquad 2 \times 17 \times 24 \times 4 \qquad 3264 & \{ \{ 18, 30, 34, 38, 50 \}, \{ 3, 23, 31, 35, 39, 47, 67 \} , \\ & \{ 16, 20, 28, 36, 44, 52, 56 \} , \{ 1, 13, 21, 25, 41, 45, 53, 65 \} \} \\
\hline 17 \qquad 2 \times 17 \times 24 \times 4 \qquad 3264 & \{ \{ 18, 30, 34, 38, 50 \}, \{ 3, 23, 27, 35, 43, 47, 67 \} , \\ & \{ 4, 16, 28, 36, 44, 56, 68 \} , \{ 13, 17, 21, 29, 37, 45, 49, 53 \} \}\\
\hline 17 \qquad 2 \times 17 \times 24 \times 4 \qquad 3264 & \{ \{ 18, 30, 34, 38, 50 \}, \{ 7, 15, 31, 35, 39, 55, 63 \} , \\ & \{ 12, 20, 24, 36, 48, 52, 60 \} , \{ 5, 9, 13, 21, 45, 53, 57,61 \} \}\\
\hline 17 \qquad  2 \times 17 \times 24 \times 4 \qquad 3264 & \{ \{ 18, 26, 30, 38, 42, 50 \}, \{ 3, 19, 31, 39, 51, 67 \} , \\ & \{ 8, 24, 32, 36, 40, 48, 64 \} , \{ 9, 13, 17, 33, 49, 53, 57 \} \}\\
\hline  \hline 19 \qquad 2 \times 19 \times 24 \times 9 \qquad 8208 & \{ \{ 18, 30, 34, 42, 46, 58 \}, \{ 15, 19, 23, 39, 55, 59, 63 \} , \\ & \{ 8, 16, 28, 36, 40, 44, 52, 64, 72 \} , \{  5, 21, 25, 29, 37, 45, 49, 53, 69 \} \} \\
\hline 19 \qquad 2 \times 19 \times 12 \times 9 \qquad 4104 &  \{ \{ 10, 30, 34, 42, 46, 66 \}, \{ 11, 31, 35, 39, 43, 47, 67 \} , \\ & \{ 4, 12, 16, 32, 40, 48, 64, 68, 76 \} , \{  1, 9, 13, 29, 37, 45, 61, 65, 73 \} \}\\
\hline 19 \qquad 2 \times 19 \times 24 \times 9 \qquad 8208 & \{
\{ 10, 30, 34, 42, 46, 66 \}, \{ 15, 27, 35, 39, 43, 51, 63 \} ,
\\  & \{ 8, 16, 32, 36, 40, 44, 48, 64, 72 \} ,
\{  5, 9, 25, 33, 37, 41, 49, 65, 69 \} \}\\
\hline  19 \qquad 2 \times 19 \times 24 \times 3 \qquad 2376 & \{ \{ 6, 10, 34, 42, 66, 70 \}, \{ 11, 23, 27, 35, 43, 51, 55, 67 \} , \\ & \{ 8, 16, 20, 36, 44, 60, 64, 72 \} ,  \{  1, 5, 9, 29, 45, 65, 69, 73 \} \} \\
\hline 19 \qquad 2 \times 19 \times 24 \times 3 \qquad 2736 & \{ \{ 6, 10, 34, 42, 66, 70 \}, \{ 7, 15, 31, 35, 43, 47, 63, 71 \} , \\ & \{ 12, 20, 24, 36, 44, 56, 60, 68 \} ,  \{  1, 5, 9, 25, 49, 65, 69, 73 \} \}\\
\hline 19 \qquad 2 \times 19 \times 24 \times 9 \qquad 8208 & \{ \{ 22, 30, 34, 42, 46, 54 \}, \{ 7, 15, 19, 35, 43, 59, 63, 71 \} , \\ & \{ 4, 16, 20, 36, 44, 60, 64, 76 \} , \{  9, 13, 17, 25, 49, 57, 61, 65 \} \}\\
\hline
\end{array}$
\end{table}

\begin{table}[h]\renewcommand{\arraystretch}{1.3}\caption{$\Z _t \times \Z_2^2$-cocyclic Hadamard matrices [21-23] (only for symmetric diagrams)}\label{sime2}\centering
$\begin{array}{c} t/ orbit/ \# H. / representative \\
\hline 21 \qquad  2 \times 21 \times 24 \times 2 \qquad 2016 \\ \{
\{  22, 26, 38, 46, 58, 62 \}, \{ 3, 15, 19, 27, 39, 47, 59, 67,
71, 83\} , \\ \{ 8, 12, 16, 24,  40, 48, 64, 72, 76, 80 \} ,
 \{ 13, 21, 25, 29, 33, 49, 53, 27, 61, 69 \} \}\\
\hline  21 \qquad 2 \times 21 \times 24 \times 6 \qquad 6048 \\ \{
\{ 10,  34, 38, 42, 46, 50, 74 \}, \{ 3, 19, 23, 31, 55, 63, 67,
83 \} , \\ \{ 12,  24, 28, 36, 44, 52, 60, 64, 76\} ,
 \{ 5, 9, 13, 21, 33, 49, 61, 69, 73, 77 \} \}\\
 \hline 21 \qquad 2 \times 21 \times 24 \times 6 \qquad 6048 \\ \{ \{ 18, 26, 38, 42, 46, 58, 66 \}, \{ 19, 23, 31, 35, 51, 55, 63, 67 \} , \\  \{ 12,  16, 32, 40, 44, 48, 56, 72, 76\} ,
 \{ 1, 5,  17, 29, 37, 45, 53, 65, 77, 81 \} \}\\
\hline 21 \qquad 2 \times 21 \times 24 \times 6 \qquad 6048 \\ \{
\{ 18, 26, 38, 42, 46, 58, 66 \}, \{ 11, 23, 27, 35, 51, 59, 63,
75 \} , \\ \{ 12,  16, 32, 40, 44, 48, 56, 72, 76\} ,
 \{ 1, 5, 9, 25, 37, 45, 57, 73, 77, 81 \} \} \\
\hline  21\qquad  2 \times 21 \times 12 \times 6 \qquad 3024 \\ \{
\{ 6, 30, 34, 38, 46, 50, 54, 78 \}, \{ 7, 31, 35, 39, 47, 51, 55,
79 \} , \\ \{ 16, 20, 36, 52, 68, 72, 80 \} ,
 \{ 13, 17, 33, 41, 49, 65, 69, 77 \} \}\\
\hline  21 \qquad 2 \times 21 \times 24 \times 6 \qquad 6048 \\ \{
\{ 18, 30, 34, 38, 46, 50, 54, 66 \}, \{ 3, 7, 31, 39, 47, 55, 79,
83 \} , \\ \{ 8, 24, 32, 36, 52, 56, 64, 80\} ,
 \{ 1, 17, 21, 29, 41, 53, 61, 65, 81 \} \}\\
 \hline  21 \qquad 2 \times 21 \times 24 \times 2 \qquad 2016 \\ \{ \{ 14, 18, 30, 38, 46, 54, 66, 70 \}, \{ 7, 15, 27, 31, 55, 59, 71, 79 \} , \\  \{ 8, 12, 20, 24, 64, 68, 72, 80 \} ,
 \{ 13, 21, 25, 37, 41, 45, 57, 61, 69 \} \}\\
\hline  \hline 23 \qquad 2 \times 23 \times 24 \times 11 \qquad
12144 \\ \{ \{  22, 30, 38, 42, 50, 54, 62, 70 \}, \{ 15, 19, 23,
35, 47, 59, 71, 75, 79\} , \\  \{ 12, 20, 28, 40, 44, 52, 56, 68,
76, 84 \} ,
 \{ 9, 13, 25, 29, 33, 57, 61, 65, 77, 81 \} \}\\
  \hline
\end{array}$
\end{table}

\begin{table}[h]\renewcommand{\arraystretch}{1.1}\caption{Full orbits from Williamson type matrices}\label{kou}\centering
$\begin{array}{c|l|r}
t & 4t=A^2+B^2+C^2+D^2 & \sharp orbit  (compl./rot./swapp./dilat.)   \\
\hline
25 & 100=9^2+3^2+3^2+1^2 & 2\times 25\times 24\times 10=12000  \\
\cline{2-3} & 100=7^2+7^2+ 1^2 + 1^2 &  2\times 25\times 12\times 5=3000  \\
&  &  2\times 25\times 24\times 5=6000  \\
& & 2\times 25\times 24\times 5=6000\\
\cline{2-3} & 100=7^2+5^2+5^2+1^2 & 2\times 25\times 24\times 10=12000 \\
& & 2\times 25\times 24\times 10=12000 \\
& & 2\times 25\times 24\times 5=6000 \\
\cline{2-3} & 100=5^2+5^2+5^2+5^2 & 2\times 25\times 24\times 5=6000 \\
& & 2\times 25\times 24\times 5=6000 \\
& & 2\times 25\times 24\times 5=6000 \\
\hline
27 & 108=9^2+5^2+1^2+1^2 & 2\times 27\times 24\times 9=11664  \\
   & & 2\times 27\times 24\times 9=11664 \\
\cline{2-3} & 108=7^2+7^2+ 3^2 + 1^2 &  2\times 27\times 12\times 9=5832  \\
&  &  2\times 27\times 24\times 9=11664  \\
& & 2\times 27\times 24\times 9=11664 \\
\cline{2-3} & 108=7^2+5^2+5^2+3^2 & 2\times 27\times 24\times 9=11664 \\
\hline
29 & 116=9^2+5^2+3^2+1^2 & 2\times 29\times 24\times 14=19488  \\
\hline
31 & 124=7^2+7^2+ 5^2 + 1^2 &  2\times 31\times 24\times 15=22320  \\
\cline{2-3} & 124=7^2+5^2+5^2+5^2 & 2\times 31\times 12\times 15=11160 \\
\hline
33 & 132=11^2+3^2+1^2+1^2 & 2\times 33\times 24\times 10=15840  \\
\cline{2-3} & 132=9^2+7^2+ 1^2 + 1^2 &  2\times 33\times 24\times 10=15840  \\
\cline{2-3} & 132=9^2+5^2+5^2+1^2 & 2\times 33\times 24\times 10=15840 \\
 & & 2\times 33\times 24\times 10=15840\\
\cline{2-3} & 132= 7^2+7^2+5^2+ 3^2 & 2\times 33\times 24\times 10=15840 \\
\hline
37 & 148=11^2+3^2+3^2+3^2 & 2\times 37\times 24\times 3=5328  \\
\cline{2-3} & 148=9^2+7^2+ 3^2 + 3^2 &  2\times 37\times 12\times 18=15984  \\
\cline{2-3} & 148=7^2+7^2+ 5^2 + 5^2 &  2\times 37\times 24\times 9=15984  \\
&  &  2\times 37\times 24\times 9=15984  \\
\hline 39 & 156=9^2+5^2+5^2+5^2& 2\times 39\times 24\times 4=7488 \\
\hline
41 & 164=9^2+9^2+1^2+1^2 & 2\times 41\times 12\times 5=4920  \\
\hline
43 & 172=7^2+7^2+7^2+5^2 & 2\times 43\times 24\times 7=14448  \\
\hline
45 & 180=9^2+7^2+5^2+5^2 & 2\times 45\times 12\times 12=12960  \\
\hline
49 & 196=9^2+9^2+5^2+3^2 & 2\times 49\times 12\times 21=24696  \\
\hline
51 & 204=11^2+9^2+ 1^2 + 1^2 &  2\times 51\times 12\times 16=19584  \\
\cline{2-3} & 204=11^2+7^2+5^2+3^2 & 2\times 51\times 24\times 16=39168 \\
\hline
55 & 220=11^2+9^2+3^2+3^2 & 2\times 55\times 12\times 20=26400  \\
\hline
57 & 228=9^2+7^2+7^2+7^2 & 2\times 57\times 12\times 18=24624  \\
\hline
61 & 244=11^2+11^2+1^2+1^2 & 2\times 61\times 12\times 15=21960  \\
\hline
63 & 252=11^2+11^2+3^2+1^2 & 2\times 63\times 12\times 6=9072  \\
\hline
\end{array}$
\end{table}


\end{article}

\end{document}